\documentclass[11pt, twoside, eqno]{article}
\usepackage{amssymb}
\usepackage{mathrsfs}
\usepackage{amsmath}
\usepackage{amsfonts}
\usepackage{amsthm}

\allowdisplaybreaks \numberwithin{equation}{section}

\def\rr{\mathbb{R}}
\def\rn{{{\rr}^n}}

\def\min{\mathop\mathrm{\,min\,}}

\def\dint{\displaystyle\int}
\def\dfrac{\displaystyle\frac}
\def\dprod{\displaystyle\prod}

\newtheorem{thm}{\noindent\rm\bf Theorem}[section]
\newtheorem{lem}{\noindent\rm\bf Lemma}[section]
\newtheorem*{rem*}{Remark}

\begin{document}

\baselineskip=15pt
\renewcommand{\arraystretch}{2}
\arraycolsep=1pt

\title{\bf \Large Weighted norm inequalities for multilinear singular integral operators and applications\footnotetext{\hspace{-0.6cm}Research of the first author
was supported by National Natural Science Foundation of China under Grant \#10971228.
Research of the second author was supported by National Science Council of Taiwan under Grant \#NSC 100-2115-M-008-002-MY3.}}
\author{\vspace{-0.2cm}\bf Guoen Hu  \\
\vspace{-0.5cm}\small\it Department of Applied Mathematics,
Zhengzhou Information Science and Technology Institute\\
\vspace{-0.5cm}\it\small P.O. Box 1001-747, Zhengzhou 450002, China\\
\it\small E-mail: guoenxx@yahoo.com.cn\\
\vspace{-0.2cm \bf Chin-Cheng Lin}\\
\vspace{-0.5cm}{\it\small Department of Mathematics, National Central University}\\
\vspace{-0.5cm}{\it\small Chung-Li 320, Taiwan}\\
{\it\small E-mail: clin@math.ncu.edu.tw}}

\date{}

\maketitle

\begin{center}
\begin{minipage}{15.2cm}
\small{\bf Abstract.}\quad
In this paper, weighted norm inequalities with $A_p$ weights are established for the multilinear singular integral operators whose kernels satisfy $L^{r'}$-H\"ormander regularity condition.
As applications, we recover a weighted estimate for the multilinear Fourier multiplier obtained by Fujita and Tomita, and obtain several new weighted estimates for the multilinear Fourier multiplier as well.

\medskip


\noindent{\bf MR(2000) Subject Classification.}\quad 42B20, 42B30.
\end{minipage}
\end{center}

\section{Introduction.}

The study of the multilinear singular integral operators was
originated by Coifman and Meyer in their celebrated work \cite{cm1,cm2}. Let $m\in \mathbb N$ and
$K(x;\,y_1,\ldots,y_m)$ be a locally integrable function defined
away from the diagonal $x=y_1=y_2=\dots=y_m$ in $(\rn)^{m+1}$. An
operator $T$, defined on $m$-fold products of $\mathscr{S}(\rn)$
(Schwartz space) and taking values in the space of tempered
distributions, is said to be an $m$-linear singular integral
operator with kernel $K$ if $T$ is $m$-linear and satisfies that
\begin{eqnarray}\label{eq1.1}
T(f_1,\dots,f_m)(x)=\int_{(\rn)^m}K(x;\,y_1,\ldots,y_m)f_1(y_1)\ldots f_m(y_m)d\vec{y},
\end{eqnarray}
for bounded functions $f_1,\ldots,f_m$ with compact supports, and
$x\in\rn\backslash \cup_{j=1}^m{\rm supp}\,f_j$, where $d\vec{y}=dy_1\ldots dy_m$. Operators of this type plays an important role in multilinear harmonic analysis. When $T$ is an
$m$-linear Calder\'on-Zygmund operator; that is, $T$ is bounded from
$L^{q_1}(\rn)\times\dots\times L^{q_m}(\rn)$ to $L^{q}(\rn)$ for some
$q_1,\dots,q_m\in [1,\,\infty]$ and $q\in (0,\,\infty)$ with
$1/q=\sum_{1\leq k\le m}1/q_k$, and the associated kernel $K$ is an
$m$-Calder\'on-Zygmund kernel, Grafakos and Torres \cite{gt1}
considered the endpoint estimate for $T$ on the space of type
$L^{1}(\rn)\times\dots\times L^{1}(\rn)$, and established a $T1$
type theorem for the operator $T$. Grafakos and Kalton \cite{gk}
proved that the multilinear Calder\'on-Zygmund operator is bounded
from the products of Hardy spaces into Lebesgue spaces. See also
\cite{gt2,h,l} for more results on the multilinear
Calder\'on-Zygmund operator.

Recently Anh and Duong \cite{at}  introduced a class of
multilinear singular integral operators whose kernels satisfy that
there exist two positive constant $r\in (1,\,\infty)$ and
$\varrho\in (0,\,1]$ such that, for any ball $B$ and $x,\,x'\in B$,
\begin{eqnarray}\label{eq1.2}
&&\bigg(\int_{S_{j_{m}}(B)}\dots\int_{S_{j_1}(B)}|K(x;\,y_1,\dots,y_m)-K(x';\,y_1,\dots,y_m)|^{r'}d\vec{y}\bigg)^{1/r'}\nonumber\\
&&\qquad\lesssim\frac{|x-x'|^{\varrho}}{|2^{j^*}B|^{m/r+\varrho/n}},
\end{eqnarray}
where $j_1,\dots,j_m$ are integers with $\max_{1\leq j\leq m}j_k>0$ and $j^*=\max_{1\leq k\leq m}j_k$.
Here and what follows, we denote by $r'$ the index conjugate to $r$; that is, $r'=r/(r-1)$ for $1\le r\le\infty$. For a ball $B$, we denote $S_N(B):=2^NB\backslash 2^{N-1}B$ for $N\in \mathbb N$, and
$S_0(B):=B$. Anh and Duong considered the weighted estimates with
multiple weights for the multilinear singular integral operators
when the associated kernels satisfy (\ref{eq1.2}), and obtained several new
weighted estimates for multilinear Fourier multiplier operators.

The purpose of this paper is to establish weighted norm
inequalities for  multilinear singular integral operators whose
kernels satisfy certain $L^{r'}$-H\"ormander condition.
Before stating our results, we first recall some notations.

A function $w$ is said to be a weight if it is nonnegative and
locally integrable. Let $M$ denote the Hardy-Littlewood maximal
operator. For $r\in (0,\,\infty)$, define $M_r$ to be the operator given by
$$M_rf(x)=\{M(|f|^r)(x)\}^{1/r}.$$ For a weight $w$, the weighted
weak $L^p(\rn)$ with respect to $w$ is defined as
$$L^{p,\,\infty}(\rn,\,w)=\{f:\, \|f\|_{L^{p,\,\infty}(\rn,\,w)}<\infty\},$$
where $\|f\|_{L^{p,\,\infty}(\rn,\,w)}^p:=\sup_{\lambda>0}\lambda^pw(\{x\in\rn:\, |f(x)|>\lambda\})$.

A weight $w$ is said to belong to the Muckenhoupt
class $A_p(\rn), 1<p<\infty,$ if
$$\sup_{B}\Big(\frac{1}{|B|}w(x)dx\Big)\Big(\frac{1}{|B|}\int_{B}w(x)^{1-p'}\,dx\Big)^{p-1}<\infty,$$
where the supremum is taken over all balls $B\subset \rn$. A weight $w$ is
said to belong to the class $A_1(\rn)$ if, for any ball $B$,
$$\frac{1}{|B|}\int_{B}w(x)\,dx\lesssim \inf_{x\in B}w(x).$$
For the properties of $A_{p}(\rn)$, we refer readers to \cite{g}.

The main result of this article is the following

\begin{thm}\label{thm1.1}
Let $T$ be an $m$-linear singular integral operator with kernel $K$
in the sense of (\ref{eq1.1}). For $x, x', y_1,\dots,y_k\in \rn$, set
$$W_0(x,\,y_1,\dots,y_m;\,x')=|K(x;\,y_1,\dots,y_k)-K(x';\, y_1,\dots,y_m)|,$$
and for $y_k'\in \rn$, $1\le k\leq m$,
$$W_k(x,\,y_1,\dots,y_m;\,y_k')=|K(x;\,y_1,\dots,y_k)-K(x;\, y_1,\dots,y_k',\dots,y_m)|.$$
Let $r_1,\dots,r_m\in [1,\,\infty)$. Suppose that
\begin{itemize}
\item[\rm (i)] for any $x\in \rn$,
\begin{eqnarray*}
\sup_{R>0}R^{\frac{n}{r_1}+\dots+\frac{n}{r_m}}\Big(\int_{A_R^x}\dots\Big(\int_{A_R^x}|K(x;y_1,\dots,y_m)|^{r_m'}dy_m\Big)^{\frac{r'_{m-1}}{r_m'}}
\dots\Big)^{\frac{r_1'}{r_2'}}dy_1\Big)^{\frac{1}{r_1'}}<\infty,
\end{eqnarray*}
where $A_R^x=\{y:\,R/2\leq |y-x|\leq 2R\};$
\item[\rm (ii)] for any ball $B$ and $x,\,x'\in B$, and any $f_1,\dots,f_m$ such that
${\rm supp}\,f_k\subset \rn\backslash 4B$ for some $1\leq k\leq m$,
\begin{eqnarray*}
\int_{(\rn)^m}|W_0(x,y_1,\dots,\,y_m,\,x')| |f_1(y_1)\dots f_m(y_m)|\,d\vec{y}\lesssim
\prod_{k=1}^m\Big(M_{r_k}f_k(x)+M_{r_k}f_k(x')\Big);
\end{eqnarray*}
\item[\rm (iii)] for each integer $k$ with $1\leq k\leq m$, and each ball $B$ with radial $R$, there exists a function ${\rm H}_{k,\,B}$,
such that for function $f_k$ with ${\rm supp}\,f_k\subset B$ and any
$x\in \rn\backslash 4B$, $y_k,\,y_k'\in B$,
\begin{eqnarray*}&&\int_{(\rn)^m}|W_k(x,y_1,\dots,y_m;\,y_k)|\prod_{l=1}^m |f_l(y_l)|\,d\vec{y}\\
&&\qquad\lesssim \int_{\rn}|f_k(y_k)|{\rm
H}_{k,\,B}(x,\,y_k,\,y_k')dy_k \prod_{1\leq l\leq m,\,l\not
=k}M_{r_l}f_l(x), \end{eqnarray*}and for any integer $j_0\geq 3$,
$$\Big(\int_{S_{j_0}(B)}|{\rm H}_k(x,\,y_k,\,y_k')|^{r'_k}dx\Big)^{1/r'_k}\lesssim
\frac{R^{\varrho}}{|2^{j_0}B|^{\frac{1}{r_k}+\varrho/n}},$$
with $\varrho$ a positive constant;
\item[\rm (iv)] $T$ is bounded from $L^{q_1}(\rn)\times\dots\times
L^{q_m}(\rn)$ to $L^{q,\,\infty}(\rn)$ for some $q_1,\dots,q_m\in
[1,\,\infty]$ and $q\in (0,\,\infty)$ with $1/q=\sum_{1\leq k\leq
m}1/q_k.$
\end{itemize}
We have the following weighted estimates for $T$.
\begin{itemize}
\item[\rm (a)] If for $1\leq k\leq m$, $p_k\in (r_k,\,\infty]$
 and the weight $w_k\in A_{p_k/r_k}(\rn)$, $p\in
 (0,\,\infty)$ such that $1/p=\sum_{1\leq k\leq m}1/p_k$, then
$$\|T(f_1,\dots,f_m)\|_{L^{p}(\rn,\,\nu_{\vec{w}})}\lesssim
\prod_{k=1}^m\|f_k\|_{L^{p_k}(\rn,\,w_k)},$$
where $\nu_{\vec{w}}=\prod_{l=1}^mw_l^{p/p_l}$
$($in case $p_k=\infty$, the terms $\|f_k\|_{L^{p_k}(\rn,\,w_k)}$ on the right hand side of the inequality above is understood to be replaced by $\|f_k\|_{L^\infty(\rn)});$
\item[\rm (b)] if for some $k$ with $1\leq k\leq m$, $p_k\in (1,\,r_k')$, and $p_l\in (r_l,\,\infty]$ for any $1\leq l\leq m$
and $l\not=k,$ $w\in \cap_{1\leq l\leq m,\,l\not
=k}A_{p_l/r_l}(\rn)$ and $w^{1-p_k'}\in A_{p_k'/r_k}(\rn)$, then
$$\|T(f_1,\dots,f_m)\|_{L^{p}(\rn,\,w)}\lesssim
\prod_{l=1}^m\|f_l\|_{L^{p_l}(\rn,\,w)};$$
\item[\rm (c)] if $p_l\in (r_l,\,\infty]$ for any $1\leq l\leq m$
and $l\not=k,$  $w^{r_k}\in A_1(\rn)$, then
$$\|T(f_1,\dots,f_m)\|_{L^{p,\,\infty}(\rn,\,w)}\lesssim\|f_k\|_{L^1(\rn,\,w)}
\prod_{1\leq l\leq m,\,l\not =k}\|f_l\|_{L^{p_l}(\rn,\,w)}$$
with $1/p=1+\sum_{1\leq l\leq m,\,l\not =k}1/p_l$.
\end{itemize}
\end{thm}

We now consider the multilinear Fourier multiplier operator. Let
$\sigma\in L^{\infty}(\mathbb{R}^{nm})$. Define the $m$-linear
Fourier multiplier operator $T_\sigma$ by
\begin{eqnarray}\label{eq1.3}
T_{\sigma}(f_1,\dots,f_m)(x)&&=\int_{(\rn)^m}{\rm exp}(2\pi ix(\xi_1+\dots+\xi_m)) \nonumber \\
&&\qquad\qquad \times\sigma(\xi_1,\dots,\xi_m)\widehat{f_1}(\xi_1)\dots\widehat{f_m}(\xi_m)d\vec{\xi}
\end{eqnarray}
for $f_1,\dots,f_m\in \mathscr{S}(\rn)$, where $``\ \widehat{}\,\ "$ denotes the Fourier transform.
Coifman and Meyer \cite{cm2} proved that if $\sigma\in
C^{s}(\mathbb{R}^{nm}\backslash\{0\})$ satisfies
$$|\partial^{\alpha_1}_{\xi_1}\dots\partial^{\alpha_m}_{\xi_m}\sigma(\xi_1,\dots,\xi_m)|\leq
C_{\alpha_1,\dots,\alpha_m}(|\xi_1|+\dots+|\xi_m|)^{-(|\alpha_1|+\dots+|\alpha_m|)}$$
for all $|\alpha_1|+\dots+|\alpha_N|\leq s$ with $s\geq 2mn+1$, then
$T_m$ is bounded from $L^{p_1}(\rn)\times\dots\times L^{p_m}(\rn)$ to
$L^p(\rn)$ for all $1<p_1,\dots,p_m,\,p<\infty $ with
$1/p=\sum_{1\leq k\leq m}1/p_k$. For the case of $s\geq nm+1$,
Grafakos-Torres \cite{gt1} and Kenig-Stein \cite{ks} (for
$m=2)$ improved Coifman and Meyer's multiplier theorem to the
indices $1/m\leq p\leq 1$ by using the multilinear Calder\'on-Zygmund
operator theory.  An important progress in this area was given by Tomita.
Let $\Phi\in \mathscr{S}(\mathbb{R}^{nm})$ satisfy
\begin{eqnarray}\label{Phi}
\left\{
\begin{array}{l}\displaystyle\text{supp}\,\Phi\subset \bigg\{(\xi_1,\dots,\xi_m):\,1/2\leq \sum_{k=1}^m|\xi_k|\leq 2\bigg\};\\
            \displaystyle\sum_{\kappa\in \mathbb{Z}}\Phi(2^{-\kappa}\xi_1,\dots,2^{-\kappa}\xi_m)=1\qquad\text{for all}\ (\xi_1,\dots,\xi_m)\in \mathbb{R}^{nm}\backslash \{0\}.
\end{array}\right.
\end{eqnarray}
Set
\begin{eqnarray}\label{eq1.5}
\sigma_{\kappa}(\xi_1,\dots,\xi_m)=\Phi(\xi_1,\dots,\xi_m)\sigma(2^{\kappa}\xi_1,\dots,
2^{\kappa}\xi_m).
\end{eqnarray}
Tomita \cite{T} proved that if
\begin{eqnarray}\label{eq1.6}
\sup_{\kappa\in\mathbb{Z}}\int_{(\rn)^m}(1+|\xi_1|^2+\dots+|\xi_m|^2)^{s}|
\widehat{\sigma}_{\kappa}(\xi_1,\dots,\xi_m)|^2d\vec{\xi}<\infty
\end{eqnarray}
for some $s>mn/2$, then $T_\sigma$ is bounded from
$L^{p_1}(\rn)\times\dots\times L^{p_m}(\rn)$ to $L^p(\rn)$
provided $p_1,\dots,p_m,\,p\in (1,\,\infty)$ and $1/p=\sum_{1\leq
k\leq m}1/p_k$. Grafakos and Si \cite{gs} considered the mapping
properties from $L^{p_1}(\rn)\times\dots\times L^{p_m}(\rn)$ to
$L^{p}(\rn)$ for $T_\sigma$ when $\sigma$ satisfies (\ref{eq1.6})
and $p\leq 1$. Miyachi and Tomita \cite{mt} considered the problem
to find minimal smoothness condition for bilinear Fourier
multiplier. Let$$\|\sigma_{\kappa}\|_{W^{s_1,\dots,
s_m}(\mathbb{R}^{nm})}=\Big(\int_{\mathbb{R}^{nm}}\langle
\xi_1\rangle^{2s_1}\dots\langle\xi_m\rangle^{2s_m}|\widehat{\sigma}_{\kappa}(\xi_1,\dots,\xi_m)|^2d\vec{\xi}\Big)^{1/2},$$
where $\langle\xi_k\rangle:=(1+|\xi_k|^2)^{1/2}$. Miyachi and Tomita
\cite{mt} proved that if
$$\sup_{\kappa\in \mathbb{Z}}\|\sigma_\kappa\|_{W^{s_1,\,s_2}(\mathbb{R}^{2n})}<\infty\qquad\text{for}\ s_1, s_2>n/2,$$
then $T_{\sigma}$ is is bounded from
$L^{p_1}(\rn)\times L^{p_2}(\rn)$ to $L^{p}(\rn)$ for any $p_1,\,p_1\in
(1,\,\infty)$ and $p\geq 2/3$ with $1/p=1/p_1+1/p_2.$ Moreover,
they also gives minimal smoothness condition for which $T_{\sigma}$
is bounded from $H^{p_1}(\rn)\times H^{p_2}(\rn)$ to $L^p(\rn)$. It
should be pointed out that the argument used in \cite{mt} applies to
the case $m>2$.  As an application of Theorem \ref{thm1.1}, we have

\begin{thm}\label{thm1.2}
Let $\sigma$ be a multiplier satisfying
\begin{eqnarray}\label{eq1.7}
\sup_{\kappa\in \mathbb{Z}}\|\sigma_\kappa\|_{W^{s_1,\dots,s_m}(\mathbb{R}^{mn})}<\infty
\end{eqnarray}
for $s_1,\dots,s_m\in (n/2,\,n]$ and $T_{\sigma}$ be the operator
defined by (\ref{eq1.3}). Set $t_k=n/s_k$. We have the following weighted
estimates for $T_{\sigma}$.
\begin{itemize}
\item[\rm (a)] If $p_k\in (t_k,\,\infty]$
 and the weight $w_k\in A_{p_k/t_k}(\rn)$ for $1\leq k\leq m$, and $p\in
 (0,\,\infty)$ such that $1/p=\sum_{1\leq k\leq m}1/p_k$, then
$$\|T_{\sigma}(f_1,\dots,f_m)\|_{L^{p}(\rn,\,\nu_{\vec{w}})}\lesssim \prod_{k=1}^m\|f_k\|_{L^{p_k}(\rn,\,w_k)};$$
\item[\rm (b)] if $p_k\in (1,\,t_k')$ for some $k$ with $1\leq k\leq m$, $p_l\in (r_l,\,\infty]$ for any $1\leq l\leq m$ and $l\not=k$, $w\in \cap_{1\leq l\leq m,\,l\not
=k}A_{p_l/t_l}(\rn)$ and $w^{1-p_k'}\in A_{p_k'/t_k}(\rn)$, then
$$\|T_{\sigma}(f_1,\dots,f_m)\|_{L^{p}(\rn,\,w)}\lesssim
\prod_{l=1}^m\|f_l\|_{L^{p_l}(\rn,\,w)};$$
\item[\rm (c)] if $p_l\in (r_l,\,\infty]$ for any $1\leq l\leq m$
with $l\not=k$ and $w^{t_k}\in A_1(\rn)$, then for $1/p=1+\sum_{1\leq l\leq
m,\,l\not =k}1/p_k$,
$$\|T_{\sigma}(f_1,\dots,f_m)\|_{L^{p,\,\infty}(\rn,\,w)}\lesssim\|f_k\|_{L^1(\rn,\,w)}
\prod_{1\leq l\leq m,\,l\not =k}\|f_l\|_{L^{p_l}(\rn,\,w)}.$$
\end{itemize}
\end{thm}

\begin{rem*}\rm The conclusion (a) in Theorem \ref{thm1.2} was proved in \cite{ft}. Here we give another simpler approach, which is of independent interest.
\end{rem*}

Throughout the article, $C$ always denotes a
positive constant that may vary from line to line but remains independent of the main variables.
We use the symbol $A\lesssim B$ to denote that there exists a
positive constant $C$ such that $A\le CB$.  For any set
$E\subset\rn$, $\chi_E$ denotes its characteristic function. We use
$B(x,\,R)$ to denote a ball centered at $x$ with radius $R$,
and denote by $B_R$ the ball $B(0,\,R)$ for simplicity. For a
ball $B\subset \rn$ and $\lambda> 0$, we use $\lambda B$
to denote the ball concentric with $B$ whose radius is
$\lambda$ times of $B$'s.

\section{Proof of Theorem \ref{thm1.1}.}\label{s2}

Let $M^{\sharp}$ be the sharp maximal operator of Fefferman and
Stein; that is, for a locally integrable function $f$,
$$M^{\sharp}f(x)=\sup_{B\ni x}\frac{1}{|B|}\int_{B}|f(y)- V_B(f)|\,dy,$$
where the supremum is taken over all balls containing $x$ and $V_B(f)$ denotes the
mean value of $f$ on ball $B$. For a fixed $\delta>0$, let
$M^{\sharp}_{\delta}$ be the operator defined by
$$M_{\delta}^{\sharp}(f)(x)=\sup_{B\ni x}\inf_{c\in \mathbb{C}}\Big(\frac{1}{|B|}\int_{B}|f(y)-c|^\delta\,dy\Big)^{1/\delta}.$$

\begin{lem}\label{lem2.1}
Let $m\geq 1$ be an integer and $T$ be an $m$-linear singular
integral operator associated with kernel $K$ in the sense of (\ref{eq1.1}).
Suppose that  $T$ satisfies the assumptions (ii) and (iv) in Theorem \ref{thm1.1}.
Let $r\in (0,\,\infty)$ such that $1/r=\sum_{1\leq k\leq m}1/r_k$, where $r_k$'s are stated as in the assumption (ii) of Theorem \ref{thm1.1}. Then for any
$\delta\in (0,\min\{q,\,r/r_1,\dots,r/r_m\})$,
$$M^{\sharp}_{\delta}(T(f_1,\dots,f_m))(x)\lesssim \prod_{k=1}^mM_{\max\{q_k,\,r_k\}}f_k(x).$$
\end{lem}

\noindent
{\it Proof}. Let $x\in \rn$ and $B$ be a ball containing $x$.
Decompose $f_k$ $(1\leq k\leq m$) as
$$f_k(y)=f_k(y)\chi_{4B}(y)+f_k(y)\chi_{\rn\backslash
4B}(y)=f_k^1(y)+f_k^{2}(y).$$ The fact that $T$ is bounded from
$L^{q_1}(\rn)\times\dots\times L^{q_m}(\rn)$ to $L^{q,\,\infty}(\rn)$,
together with the argument used in the proof of the Kolmogorov
inequality, yields
\begin{eqnarray}\label{eq2.1}
\Big(\frac{1}{|B|}\int_{B}|T(f_1^1,\dots,f_m^1)(y)|^{\delta}dy\Big)^{1/\delta}
&\lesssim& \prod_{k=1}^m\Big(\frac{1}{|B|}\int_{4B}|f_k(y_k)|^{q_k}dy_k\Big)^{1/q_k} \nonumber\\
&\lesssim& \prod_{k=1}^mM_{q_k}f_k(x).
\end{eqnarray}
On the other hand, if
$i_1,\dots,i_m\in\{1,\,2\}$ and $i_k=2$ for some $k$ with $1\leq k\leq m$, we then by assumption (ii) that
\begin{eqnarray}\label{eq2.2}
|T(f_1^{i_1},\dots,f_m^{i_m})(y)-T(f_1^{i_1},\dots,f_{m}^{i_m})(y')|\lesssim
\prod_{k=1}^m\Big(M_{r_k}f_k(y)+M_{r_k}f_k(y')\Big),
\end{eqnarray}
for any $y'\in B$ such that $|
T(f_1^{i_1},\dots,f_{m}^{i_m}(y')|<\infty$. It follows from \cite{ac} that
$$M(M_{r_k}f_k)(z)\lesssim M_{r_k}f_k(z)\qquad\text{for}\ z\in\rn.$$
Both estimates (\ref{eq2.1}) and (\ref{eq2.2}) lead to
\begin{eqnarray*}&&\inf_{c\in
\mathbb{C}}\Big(\frac{1}{|B|}\int_B|T(f_1,\dots,f_m)(y)-c|^{\delta}dy\Big)^{\delta}\\
&&\qquad\lesssim
\Big(\frac{1}{|B|^2}\int_B\int_B|T(f_1,\dots,f_m)(y)-T(f_1,\dots,f_m)(y')|^{\delta}dydy'\Big)^{\delta}
\\
&&\qquad\lesssim\Big(\frac{1}{|B|}\int_B|T(f_1^1,\dots,f_m^1)(y)|^{\delta}dy\Big)^{\delta}\\
&&\qquad\quad+\sum_{i_1,\dots,i_m}^*
\Big(\frac{1}{|B|^2}\int_B\int_B|T(f_1^{i_1},\dots,f_m^{i_m})(y)-T(f_1^{i_1},\dots,f_m^{i_m})(y')|^{\delta}dydy'\Big)^{\delta}
\\
&&\qquad\lesssim\prod_{k=1}^mM_{q_k}f_k(x)+\Big(\frac{1}{|B|}\int_B\prod_{k=1}^m\{M_{r_k}f_k(y)\}^{\delta}dy\Big)^{1/\delta}
\\
&&\qquad\lesssim\prod_{k=1}^mM_{q_k}f_k(x)+\prod_{k=1}^mM_{r_k\delta/r}(M_{r_k}f_k)(x)\\
&&\qquad\lesssim\prod_{k=1}^mM_{\max\{q_k,\,r_k\}}f_k(x),\end{eqnarray*}
where, for each term in the summation $\sum_{i_1,\dots,i_m}^*$, the
set of indices $\{i_1,\dots,i_m\}\subset\{1,\,2\}$ and at least one
$i_k=2$ ($1\leq k\leq m$). This finishes the proof. \qed

\begin{lem}\label{lem2.2}
Let $p\in (0,\,\infty)$. If there exists $p_0\in (0, p)$ such that $\|M_{\delta}f\|_{L^{p_0,\,\infty}(\rn)}<\infty$, then
\begin{eqnarray}\label{eq2.3}
\|M_{\delta}f\|_{L^p(\rn)}\lesssim \|M_{\delta}^{\sharp}f\|_{L^{p}(\rn)}.
\end{eqnarray}
\end{lem}

\noindent
{\it Proof}. Note that, for $\delta\in (0,\,1)$,
$$\{M^{\sharp}(|h|^{\delta})(x)\}^{1/\delta}\lesssim M_{\delta}^{\sharp}h(x).$$
If $\|M_{\delta}f\|_{L^{p}(\rn)}<\infty$ holds, then (\ref{eq2.3})
follows from \cite[Theorem 7.4.5 and Corollary 7.4.6]{g}. On the
other hand, for any positive real number $N$,
$$\int_0^N \lambda^{p-1} |\{x\in \mathbb R^n: M_\delta f(x)> \lambda\}|\, d \delta
  \lesssim N^{p-p_0} \sup_{\lambda >0} \lambda^{p_0} |\{x\in \mathbb R^n: M_\delta f(x)> \lambda\}|.$$
Thus, by the same argument in the proof of \cite[Theorem 7.4.5]{g}, the lemma follows.
\qed

\begin{lem}\label{lem2.3}
Let $m\in \mathbb N$ and $T$ be an $m$-linear singular integral
operator associated with kernel $K$ in the sense of (\ref{eq1.1}). Suppose
that
\begin{itemize}
\item[\rm (1)] $T$ satisfies assumption (iii) in Theorem \ref{thm1.1};
\item[\rm (2)]  $T$ is bounded from
$L^{u_1}(\rn)\times\dots\times L^{u_m}(\rn)$ to
$L^{u,\,\infty}(\rn)$ for some $u_1,\dots,u_m$ with $u_k\in
[r_k,\,\infty]$ $(1\leq k\leq m)$, and $u\in (0,\,\infty)$ with
$1/u=\sum_{1\leq k\leq m}1/u_k.$
\end{itemize}
Then, for $p_1,\dots,p_m$ such that $p_k\in [r_k,\,u_k]$, $T$ is
bounded from $L^{p_1}(\rn)\times\dots\times L^{p_m}(\rn)$ to
$L^{p,\,\infty}(\rn)$ with $1/p=\sum_{1\leq k\leq m}1/p_k.$
\end{lem}

\noindent
{\it Proof}. Let $p_k\in [r_k,\,u_k]$ ($1\leq k\leq m)$ and
$$\|f_1\|_{L^{p_1}(\rn)}=\dots=\|f_m\|_{L^{p_m}(\rn)}=1.$$
Our goal is to prove that there exists a constant depending only on $n$, $m$ and $p$ such that
\begin{eqnarray}\label{eq2.5}
\big|\big\{x\in\rn:\, |T(f_1,\dots,f_m)(x)|>(C+1)\lambda\big\}\big|\lesssim
\lambda^{-p}\qquad\text{for all}\ \lambda>0.
\end{eqnarray}
To do this, we apply the Calder\'on-Zygmund decomposition. Given $\lambda>0$, applying the Calder\'on-Zygmund decomposition to $|f_k|^{p_k}$ at level $\lambda^p$, we obtain a sequence of cubes
$\{Q_{k}^j\}_j$ satisfying
$$\lambda^p<\frac{1}{|Q_k^j|}\int_{Q_k^j}|f_k(y)|^{p_k}dy\leq 2^n\lambda^p$$
and
$$|f_k(x)|\leq \lambda^{p/p_k}\qquad{\rm a.e.}\ x\in \rn\backslash \cup_{j}Q_{k}^j.$$
Set
$$\aligned
f_k^1(y)&:=f_k(y)\chi_{\rn\backslash
\cup_jQ_k^j}(y)+\sum_{j}V_{Q_k^j}(f_k)\chi_{Q_k^j}(y),\\
f^2_k(y)&:= f_k(y)-f^1_k(y)=\sum_{j}b_k^j(y),
\endaligned$$
where $b_k^j(y)=(f_k(y)-V_{Q_k^j}(f_k))\chi_{Q_k^j}(y)$.
It is well known that $f^1_k\in L^{u_k}(\rn)$ and
$$\aligned
\|f^1_k\|_{L^{u_k}(\rn)}&\lesssim \lambda^{p/p_k(1-p_k/u_k)}\|f_k\|_{L^{p_k}(\rn)}^{p_k/u_k}\lesssim \lambda^{p/p_k-p/u_k},\\
\|f^2_k\|_{L^{p_k}(\rn)}&\lesssim
\|f_k\|_{L^{p_k}(\rn)}\lesssim 1.
\endaligned$$
Recall that $T$ is bounded from $L^{u_1}(\rn)\times\dots
L^{u_m}(\rn)$ to $L^{u,\,\infty}(\rn)$. It follows that
$$|\{x\in \rn:\, |T(f^1_1,\dots,f^1_m)(x)|>\lambda\}|\lesssim
\lambda^{-u}\prod_{k=1}^m\|f^1_k\|_{L^{u_k}(\rn)}^u\lesssim
\lambda^{-p}.$$ The proof of (\ref{eq2.5}) is now reduced to proving
\begin{eqnarray}\label{eq2.6}
|\{x\in \rn:\,
\sum_{i_1,\dots,i_m}^*|T(f_1^{i_1},\dots,f_m^{i_m})(x)|>C\lambda\}|\lesssim
\lambda^{-p},
\end{eqnarray}
where, for each term $T(f_1^{i_1},\dots,f_m^{i_m})$ in the sum
$\sum_{i_1,\dots,i_m}^*$, each one of $i_1,\dots,i_m$ is either 1 or 2, and at
least one $i_k=2$.

To prove (\ref{eq2.6}), without loss of generality we may assume
$i_1=2$. Let $\Omega=\cup_{k=1}^m\cup_j8B_k^j$, where $B_k^j$ is the
smallest ball containing $Q_k^j$. For each $x\in \rn\backslash
\Omega$ and each fixed $j$, applying the vanishing moment of $b_1^j$
and the H\"older inequality, we write
\begin{eqnarray*}&&|T(b^j_1,f_2^{i_2},\dots,f_{i_m}^{i_m})(x)|\\
&&\qquad\lesssim
\int_{(\rn)^m}|K(x;y_1,\dots,y_m)-K(x;y_1^j,y_2,\dots,y_m)||b_1^j(y_1)|\prod_{k=2}^m|f_k^{i_k}(y_k)|d\vec{y}\\
&&\qquad\lesssim\int_{B_1^j}|b_1^j(y_1)|{\rm
H}_{1,\,B_1^j}(x,\,y_1,\,y_1^j)dy_1\prod_{k=2}^mM_{r_k}f_k^{i_k}(x),
\end{eqnarray*}
where $y_1^j$ is the center of $B_1^j$. Let
\begin{eqnarray}\label{eq2.7}
{\rm I}(x)=\sum_j\int_{B_1^j}|b_1^j(y_1)|{\rm H}_{1,\,B_1^j}(x,\,y_1,\,y_1^j)dy_1.
\end{eqnarray}
We then have
\begin{eqnarray}\label{eq2.8}
|T(f_1^2,f_2^{i_2},\dots,\,f_m^{i_m}(x)|\lesssim {\rm I}(x)\prod_{k=2}^m M_{r_k}f_k^{i_k}(x).
\end{eqnarray}
Recall that, for any $k$ with $2\leq k\leq m$,
$$M_{r_k}f_k^{i_k}(x)\leq M_{r_k}f_k(x)+D\lambda^{p/p_k}$$
for some constant $D$ depending only on $n$. Thus,
\begin{eqnarray}\label{eq2.9}
|\{x\in \rn:\, Mf_k^{i_k}(x)>(D+1)\lambda^{p/p_k}\}|\lesssim
\lambda^{-p}\int_{\rn}|f_k(y)|^{p_k}\,dy\lesssim\lambda^{-p}.
\end{eqnarray}
On the other hand, a trivial computation yields
\begin{eqnarray*}
\int_{\rn\backslash \Omega}{\rm I}(x)\,dx &=&\sum_j
\sum_{l=3}^{\infty}\int_{S_l(B_1^j)}\int_{B_1^j}\int_{S_l(B_1^j)}{\rm H}_{1,\,B_1^j}(x,\,y_1,\,y_1')dx|b_1^j(y_1)|dy_1\\
&\lesssim&\sum_j
\sum_{l=3}^{\infty}\int_{B_1^j}\Big(\int_{S_l(B_1^j)}|{\rm
H}_{1,\,B_1^j}(x;y_1,y_1^j)|^{r_1'}dx\Big)^{1/r_1'}|b_1^j(y_1)|dy_1
|S_l(B_1^j)|^{1/r_1}\\
&\lesssim &\sum_j\int_{B_1^j}|b_i^j(y_1)|dy_1
\end{eqnarray*}
that implies
\begin{eqnarray}\label{eq2.10}
\big|\{x\in \rn\backslash \Omega:\,{\rm I}(x)>\lambda^{p/p_1}\}\big| \leq
\lambda^{-p/p_1}\int_{\rn\backslash \Omega}{\rm I}(x)\,dx\lesssim \lambda^{-p}.
\end{eqnarray}
Combining inequalities $(\ref{eq2.8})-(\ref{eq2.10})$, we obtain (\ref{eq2.6}).
This completes the proof of Lemma \ref{lem2.3}.
\qed
\medskip

We now are ready to show the main theorem.

\medskip\noindent
{\it Proof of Theorem \ref{thm1.1}}. Since $T$ is bounded from
$L^{q_1}(\rn)\times\dots\times L^{q_m}(\rn)$ to
$L^{q,\,\infty}(\rn)$, it follows from Lemmas \ref{lem2.1} and
\ref{lem2.2} that $T$ is bounded from $L^{p_1}\times\dots\times
L^{p_m}(\rn)$ to $L^{p}(\rn)$ with $1/p=\sum_{1\leq k\leq m}1/p_k$
provided $p_k> \max\{q_k,\,r_k\}$ for $k=1,\dots,m$. Thus, by Lemma
\ref{lem2.3}, $T$ is bounded from $L^{r_1}(\rn)\times\dots \times
L^{r_m}(\rn)$ to $L^{r,\infty}(\rn)$ with $r\in (0,\infty)$ and
$1/r=\sum_{1\leq k\leq m}1/r_k$. Hence, for $\delta\in
(0,\min\{r/r_1,\dots,r/r_m\})$,
\begin{eqnarray}\label{eq2.11}
M_{\delta}^{\sharp}(T(f_1,\dots,f_m))(x)\lesssim \prod_{k=1}^{m}M_{r_k}f_k(x).
\end{eqnarray}
Applying Lemma \ref{lem2.2} again, we obtain the boundedness of $T$ from $L^{p_1}(\rn)\times\dots\times
L^{p_m}(\rn)$ to $L^{p}(\rn)$ provided $p_k\in (r_k,\infty], 1\leq k\leq m,$ and $p\in (0,\infty)$ such that $1/p=\sum_{1\leq
k\leq m}1/p_k.$

To prove conclusion (a), for $p_k\in (r_k,\infty]$ and $w_k\in A_{p_k/r_k}(\rn)$, $k=1,\dots,m$,
we claim that if $\delta$ is small enough, then for
bounded functions $f_1,\dots,f_m$ with compact supports,
\begin{eqnarray}\label{eq2.12}
\|M_{\delta}(T(f_1,\dots,f_m))\|_{L^{p}(\rn,\,\nu_{\vec{w}})}<\infty.
\end{eqnarray}
Once we prove the claim, conclusion (a) follows immediately from the
inequalities $(\ref{eq2.11})-(\ref{eq2.12})$ and the well known inequality of C\'ordoba
and Fefferman \cite{cf}.

The proof of (\ref{eq2.12}) is fairly standard. We note that
$\nu_{\vec{w}}\in A_{p/\delta}(\rn)$ for $\delta$ small enough. If
we take $R$ large enough such that $\cup_{k=1}^m{\rm
supp}\,f_k\subset B_{R}$, then
\begin{eqnarray*}\|M_{\delta}(T(f_1,\dots,f_m))\|_{L^{p}(\rn,\,\nu_{\vec{w}})}^p&\lesssim &
\|T(f_1,\dots,f_m)\|_{L^p(\rn,\,\nu_{\vec{w}})}^p\\
&=&\int_{B_{2R}}|T(f_1,\dots,f_m)(x)|^p\nu_{\vec{w}}(x)\,dx\\
&&+\int_{\rn\backslash
B_{2R}}|T(f_1,\dots,f_m)(x)|^p\nu_{\vec{w}}(x)\,dx.
\end{eqnarray*}
It is obvious that
$\int_{B_{2R}}|T(f_1,\dots,f_m)(x)|^p\nu_{\vec{w}}(x)\,dx<\infty$.
On the other hand, the size condition (i) shows that, for
$x\in\rn\backslash B_{2R}$,
\begin{eqnarray*}
|T(f_1,\dots,f_m)(x)|
&\lesssim& \int_{\{\mathbb{R}^{mn}:\frac{|x|}{2}\leq |x-y_k|\leq 2|x|,k=1,\dots,m\}}|K(x;y_1,\dots,y_m)||f_1(y_1)\dots f_m(y_m)|d\vec{y}\\
&\lesssim &\prod_{k=1}^mM_{r_k}f_k(x).
\end{eqnarray*}
Then, H\"older's inequality and the weighted boundedness of the
Hardy-Littlewood maximal operator yield
$$\int_{\rn\backslash B_{2R}}|T(f_1,\dots,f_m)(x)|^p\nu_{\vec{w}}(x)\,dx\lesssim
\prod_{k=1}^m\|f_k\|_{L^{p_k}(\rn,\,w_k)}^p.$$ This leads to our
claim directly.

For conclusion (b), we consider the case $k=1$ only. Let $1<p_1<r_1'$, $p_l> r_l$ for $2\leq l\leq m$,
$w\in\cap_{l=2}^mA_{p_l/r_l}(\rn)$ and $w^{1-p_1'}\in
A_{p_1'/r_1}(\rn)$. Choose points
$(1/p_1^1,\dots,1/p_m^1,1/p^1),\dots,\break (1/p_1^{m+1},\dots,1/p_m^{m+1},1/p^{m+1})$,
such that $1/p^j=\sum_{1\leq k\leq m}1/p_{k}^j$ for any $1\leq j\leq m+1$, and
\begin{itemize}
\item[\rm (1)] $(1/p_1,\dots,1/p_m,\,1/p)$ is in the open convex hull of the points $(1/p_1^1,\dots,1/p_m^1,\,1/p^1),$ \dots, $(1/p_1^{m+1},\dots,1/p_m^{m+1},1/p^{m+1})$;
\item[\rm (2)] for each $i\in \{1,\dots, m+1\}$, either $p_l^i> r_l$ for all $1\leq l\leq m$ and
$w\in \cap_{l=1}^m A_{p_l^i/r_l}(\rn)$, or $1<p_1^i<r_1'$,
$p_l^i>r_l$ for $2\leq l\leq m$, $w^{1-(p_1^i)'}\in
A_{(p_1^i)'/r_1}(\rn)$ and $w\in \cap_{l=2}^m A_{p_l^i/r_l}(\rn)$.
\end{itemize}
Thus, by the multilinear Marcinkiewicz interpolation theorem (cf. \cite{g}), it suffices to prove the boundedness of $T$ from $L^{p_1}(\rn,\,w) \times\dots\times L^{p_m}(\rn,\,w)$ to $L^{p,\infty}(\rn,\,w)$ whenever $p_1\in (1,\,r_1')$, $p_l> r_l$ for $2\leq l\leq m$, $w^{1-p_1'}\in A_{p_1'/r_1}(\rn)$ and $w\in A_{p_k/r_k}(\rn)$ for $2\leq k\leq m$.

For $f_k\in L^{p_k}(\rn), 1\leq k\leq m$, with
$$\|f_1\|_{L^{p_1}(\rn,\,w)}=\dots=\|f_m\|_{L^{p_m}(\rn,\,w)}=1,$$
applying the weighted Calder\'on-Zygmund decomposition to
$|f_1|^{p_1}$ at level $\lambda^p$, we obtain a sequence of cubes
$\{Q_1^j\}_{j}$ such that
$$\lambda^p\leq \frac{1}{w(Q_1^j)}\int_{Q_1^j}|f_1(y)|^{p_1}w(y)\,dy\leq 2^n\lambda^{p},$$
$$|f_1(x)|\leq \lambda^{p/p_1}\qquad\text{a.e.}\ x\in \rn\backslash \cup_{j}Q_1^j.$$
Let $f_1^1, f_1^2$ be the functions given in the proof of Lemma
\ref{lem2.3} and $\Omega=\cup_{j}8B_1^j$, where $B_1^j$ is the ball
circumscribed on $Q_1^j$. Then $w(\Omega)\lesssim \lambda^{-p}$.
Since $w^{1-p_1'}\in A_{p_1'/r_1}(\rn)$, we can choose $t_1$ large
enough such that $w\in A_{t_1/r_1}(\rn)$. Thus, by conclusion (a),
\begin{eqnarray}\label{eq2.13}
w(\{x\in \rn:\, |T(f_1^1,\,f_2,\dots,f_m)(x)|>\lambda\})
&\lesssim& \lambda^{-t}\prod_{k=2}^m\|f_k^1\|_{L^{p_k}(\rn,\,w)}^t\|f_1^1\|_{L^{t_1}(\rn,\,w)}^t\nonumber\\
&\lesssim& \lambda^{-p},
\end{eqnarray}
where $1/t=1/t_1+\sum_{2\leq k\leq m}1/p_k$.

To estimate $T(f_1^2,f_2,\dots,f_m)(x)$, we employ the idea used in \cite{kw}. Similar to the proof of Lemma \ref{lem2.3}, for $x\in\rn\backslash \Omega$,
$$|T(f_1^2,f_2^,\dots,f_m)(x)|\lesssim\sum_j\int_{B_1^j}{\rm H}_{1,\,B_1^j}(x;y_1,y_1^j)|b_1^j(y_1)|dy_1\prod_{k=2}^m M_{r_k}f_k(x).$$
By the weighted boundedness of the Hardy-Littlewood maximal operator, if $w\in A_{p_k/r_k}(\rn)$, then
\begin{eqnarray}\label{eq2.14}
w(\{x\in \rn:\, M_{r_k}f_k(x)>(D+1)\lambda^{p/p_k}\})\lesssim \lambda^{-p}\int_{\rn}|f_k(y)|^{p_k}w(y)\,dy.
\end{eqnarray}
Let ${\rm I}(x)$ be given in (\ref{eq2.7}). A duality argument shows
that if $w^{1-p_1'}\in A_{p_1'/r_1}(\rn)$ and $F\in
L^{p_1'}(\rn\backslash \Omega,\,w)$ with
$\|F\|_{L^{p_1'}(\rn\backslash \Omega,\,w^{1-p_1'})}\leq 1$, then
\begin{eqnarray*}
\Big|\int_{\rn\backslash \Omega}{\rm I}(x)F(x)\,dx\Big|
&\lesssim&\sum_j\int_{\rn\backslash \Omega}\int_{B_1^j}{\rm H}_{1,\,B_1^j}(x;y_1,y_1^j)|b_1^j(y_1)|dy_1F(x)dx\\
&\lesssim&\sum_j \sum_{l=3}^{\infty}\int_{B_1^j}\Big(\int_{S_l(B_1^j)}|{\rm
H}_{1,\,B_1^j}(x;y_1,y_1^j)|^{r_1'}dx\Big)^{1/r_1'}|b_1^j(y_1)|dy_1 \\
&&\qquad\times\inf_{y\in B_1^j}M_{r_1}F(y) |S_l(B_1^j)|^{1/r_1}\\
&\lesssim&\int_{\rn}|f_1(y)|M_{r_1}F(y)\,dy\\
&\lesssim& 1.
\end{eqnarray*}
Therefore,
\begin{eqnarray}\label{eq2.15}
w(\{x\in\rn\backslash \Omega:\,{\rm I}(x)>\lambda^{p/p_1})\lesssim
\lambda^{-p}\int_{\rn\backslash \Omega}{\rm I}(x)^{p_1}w(x)\,dx\lesssim \lambda^{-p}.
\end{eqnarray}
Estimates $(\ref{eq2.13})-(\ref{eq2.15})$ give us
$$w(\{x\in \rn:\, |T(f_1,\dots,f_m)|>C\lambda\})\lesssim \lambda^{-p}.$$

We turn the attention to conclusion (c), and consider the case $k=1$ only. Recall that $w^{r_1}\in A_ 1(\rn)$ implies that $w\in A_{r}(\rn)$ for any $r\geq 1$. A computation shows that
\begin{eqnarray*}
\int_{\rn\backslash \Omega}{\rm I}(x)w(x)\,dx
&\lesssim& \sum_j \sum_{l=3}^{\infty}\int_{B_1^j}\Big(\int_{S_l(B_1^j)}|{\rm
H}_{1,\,B_1^j}(x;y_1,y_1^j)|^{r_1'}dx\Big)^{1/r_1'}|b_1^j(y_1)|dy_1 \\
&&\qquad\times\Big(\int_{S_l(B_1^j)}w(x)^{r_1}\,dx\Big)^{1/r_1}\\
&\lesssim& \sum_j\int_{B_1^j}|b_i^j(y_1)|w(y)dy_1,
\end{eqnarray*}
which implies conclusion (c) by applying the Calder\'on-Zygmund decomposition and the estimates used in the proof of conclusion (b).
\qed

\section{Proof of Theorem \ref{thm1.2}.}\label{s3}

We start with several preliminary lemmas.

\begin{lem}\label{lem3.1}
Let $\sigma_{\kappa}$ be defined in (\ref{eq1.5}), $q_1,\dots,q_m\in
[2,\infty)$, and $s_1,\dots,s_m\geq 0$. Then
\begin{eqnarray*}&&\Big(\int_{\rn}\dots\Big(\int_{\rn}|\widehat{\sigma}_{\kappa}(\xi_1,\dots\xi_m)|^{q_1}\langle
\xi_1\rangle^{s_1}d\xi_1\Big)^{q_2/q_1}\langle\xi_2\rangle^{s_2}d\xi_2\Big)^{q_3/q_2}\dots\langle\xi_m\rangle^{s_m}d\xi_m\Big)^{1/q_m}\\
&&\qquad\lesssim\|\sigma_{\kappa}\|_{W^{s_1/q_1,\dots,s_m/q_m}(\mathbb{R}^{mn})}.
\end{eqnarray*}\end{lem}

For the proof of Lemma \ref{lem3.1}, see Appendix A in \cite{ft}.

\begin{lem}\label{lem3.2}
Let $s_1,\dots,s_m\in\mathbb{R}$, and $\alpha_1\dots,\alpha_m\in\mathbb{Z}_+^n$ be multi-indices. Set
$$\zeta_{\kappa}^{\alpha_1,\dots,\alpha_m}(\xi_1,\dots,\xi_m):=\xi_1^{\alpha_1}\dots\xi_m^{\alpha_m}\sigma_{\kappa}(\xi_1,\dots,\xi_m).$$
Then
$$\|\zeta_{\kappa}^{\alpha_1,\dots,\alpha_m}\|_{W^{s_1,\dots,s_m}(\mathbb{R}^{mn})}\lesssim
\sup_{l\in\mathbb{Z}}\|\sigma_{l}\|_{W^{s_1,\dots,s_m}(\mathbb{R}^{mn})}.$$
\end{lem}

This lemma was given in \cite[Remark 2.5]{mt}.

Let $\sigma\in L^{\infty}(\mathbb{R}^{mn})$ and $\Phi\in \mathscr{S}(\mathbb{R}^{mn})$ satisfy (\ref{Phi}). Define
$$\widetilde{\sigma}_{\kappa}(\xi_1,\dots,\xi_m)=\Phi(2^{-\kappa}\xi_1,\dots,2^{-\kappa}\xi_m)\sigma(\xi_1,\dots,\xi_m).$$
Then
$$\widetilde{\sigma}_{\kappa}(\xi_1,\dots,\xi_m)=\sigma_{\kappa}(2^{-\kappa}\xi_1,\dots,2^{-\kappa}\xi_m)$$
and
$${\cal F}^{-1}\widetilde{\sigma}_{\kappa}(\xi_1,\dots,\xi_m)=2^{\kappa nm}{\cal F}^{-1}\sigma_{\kappa}(2^{\kappa}\xi_1,\dots,2^{\kappa}\xi_m),$$
where ${\cal F}^{-1}$ denotes the inverse Fourier transform. For an integer $k$ with
$1\leq k\leq m$ and $x,\,y_1,\dots,y_m,\,y_k',\,x'\in \rn$, let
$$W_{0,\,\kappa}(x,\,y_1,\dots,y_m;\, x')={\cal F}^{-1}\widetilde{\sigma}_{\kappa}(x-y_1,\dots,x-y_m)-{\cal
F}^{-1}\widetilde{\sigma}_{\kappa}(x'-y_1,\dots,x'-y_m),$$
$$W_{k,\,\kappa}(x,\,y_1,\dots,y_m;\, y_k')={\cal F}^{-1}\widetilde{\sigma}_{\kappa}(x-y_1,\dots,x-y_m)-
{\cal
F}^{-1}\widetilde{\sigma}_{\kappa}(x-y_1,\dots,x-y_k',\dots,x-y_m).$$
Also, for $r_1,\dots,r_m\in (1,\infty)$, write
\begin{eqnarray*}
A_{k,\kappa}&=&\Big(\dint_{S_{j_0}(B)}\Big(\dint_{E_{j_1}^R(x)}\dots\Big(\dint_{E_{j_{k-1}}^R(x)}
\Big(\dint_{E_{j_{k+1}}^R(x)} \dots\Big(\int_{E_{j_{m}}^R(x)} \\
&&\quad|W_{k,\,\kappa}
(x,y_1,\dots,y_m;y_k')|^{r_m'}dy_m\Big)^{\frac{r_{m-1}'}{r_m'}}\dots dy_{k+1}\Big)^{\frac{r'_{k}}{r'_{k+1}}}
dy_{k-1}\Big)^{\frac{r'_{k-2}}{r'_{k-1}}}
\dots dy_1\Big)^{\frac{r_k'}{r_1'}}dx\Big)^{\frac{1}{r_k'}},
\end{eqnarray*}
where $E_{0}^R(x)=B(x,\,R)$ and $E_j^R(x)=2^{j}B(x,\,R)\backslash
2^{j-1}B(x,\,R)$ for $j\in \mathbb N$.

\begin{lem}\label{lem3.3}
Let $m$ and $k$ be positive integers with $k\leq m$, $\sigma$ be a
multiplier satisfying (\ref{eq1.7}) for some $s_1,\dots,s_m\in (n/2,\,n]$, and $r_1,\dots,r_m\in (1,\,2]$.
Then, for any ball $B$ with radial $R$, $y_k,\,y_k'\in\frac{1}{4}B$, $j_0\in \mathbb N$ and
nonnegative integers $j_1,\dots,j_{k-1},\,j_{k+1},\dots,j_m$,
\begin{eqnarray}\label{eq3.1}
A_{k,\,\kappa}\lesssim
\dfrac{R2^{-\kappa(s_1+\dots+s_m-n/r_1-\dots-n/r_m-1)}}{|2^{j_0}B|^{s_k/n}\prod_{1\leq
i\leq m,\,i\not =k}(2^{j_i}R)^{s_i}}\qquad\text{if}\ \
2^{\kappa}R<1.
\end{eqnarray}
\end{lem}

\noindent
{\it Proof}. We prove (\ref{eq3.1}) for $k=1$ only. Let $B_R=B(0,\,R)$.
Step 1: we first consider the case that all of $j_0,\,j_2,\dots,j_m$ are
positive. Write
\begin{eqnarray*}
&&|{\cal F}^{-1}\widetilde{\sigma}_\kappa(z_1,z_2,\dots,z_m)-{\cal
F}^{-1}\widetilde{\sigma}_\kappa(z_1+y_1-y_1',z_2,\dots,z_m)|\\
&&\qquad=2^{\kappa nm}\left|{\cal
F}^{-1}{\sigma}_\kappa(2^{\kappa}z_1,\dots,2^{\kappa}z_m)-{\cal
F}^{-1}{\sigma}_\kappa(2^{\kappa}z_1+2^{\kappa}(y_1-y_1'),\,2^{\kappa}z_2,\dots,2^{\kappa}z_m)\right|\\
&&\qquad\leq 2^{\kappa
nm}\sum_{|\alpha|=1}|2^{\kappa}(y_1-y_1')|^{\alpha}\int_{0}^1|\partial^{\alpha,\,0,\dots,0}{\cal
F}^{-1}{\sigma}_{\kappa} (
2^{\kappa}(z_1+\theta(y_1-y_1')),\,2^{\kappa}z_2\dots,2^{\kappa}\,z_m)|d\theta.
\end{eqnarray*}
Let
$$\phi_{\kappa,\,y_1-y_1'}(\theta;\,z_1,\dots,z_m)=\partial^{\alpha,\,0,\dots,0}{\cal
F}^{-1}{\sigma}_{\kappa} (
2^{\kappa}(z_1+\theta(y_1-y_1')),\,2^{\kappa}z_2\dots,2^{\kappa}\,z_m)|.$$
By the Minkowsky inequality, Lemmas \ref{lem3.1} and \ref{lem3.2},
\begin{eqnarray*}
A_{1,\kappa}
&\lesssim&\sum_{|\alpha|=1}\Big(\int_{C_{j_0}}\Big(\int_{S_{j_2}(B_R)}\dots\Big(\int_{S_{j_m}(B_R)}\Big(\int^1_0
|\phi_{\kappa,\,y_1-y_1'}(\theta;\,z_1,\dots,z_m)|d\theta\Big)^{r_m'}dz_m\Big)^{\frac{r_{m-1}'}{r_m'}}\\
&&\qquad\dots\Big)^{\frac{r_1'}{r_2'}}dz_1\Big)^{\frac{1}{r_1'}}2^{\kappa
nm}2^{\kappa}R\\
&\lesssim&\sum_{|\alpha|=1}\int_0^1\Big(\int_{C_{j_0}}\Big(\int_{S_{j_2}(B_R)}\dots\Big(\int_{S_{j_m}(B_R)}
\Big|\phi_{\kappa,\,y_1-y_1'}(\theta;\,z_1,\dots,z_m)\Big|^{r_m'}dz_m\Big)^{\frac{r_{m-1}'}{r_m'}}\\
&&\qquad \dots\Big)^{\frac{r_1'}{r_2'}}dz_1\Big)^{\frac{1}{r_1'}}d\theta2^{\kappa nm}2^{\kappa}R \\
&\lesssim&\sum_{|\alpha|=1}\Big(\int_{C_{j_0}}\Big(\int_{S_{j_2}(B_R)}\dots\Big(\int_{S_{j_m}(B_R)}
\Big|\partial^{\alpha,\,0,\dots,0}{\cal F}^{-1}{\sigma}_{\kappa} (
2^{\kappa}z_1,\dots,2^{\kappa}z_m)\Big|^{r_m'}dz_m\Big)^{\frac{r_{m-1}'}{r_{m}'}}\\
&&\qquad \dots\Big)^{\frac{r_1'}{r_2'}}dz_1\Big)^{\frac{1}{r_1'}}2^{\kappa nm}2^{\kappa}R\\
&\lesssim
&\sum_{|\alpha|=1}\Big(\int_{C_{j_0}}\Big(\int_{S_{j_2}(B_R)}\dots\Big(\int_{S_{j_m}(B_R)}
\Big|\partial^{\alpha,\,0,\dots,0}{\cal F}^{-1}{\sigma}_{\kappa} (
2^{\kappa}z_1,\dots,2^{\kappa}z_m)\Big|^{r_m'} |z_m|^{r_m's_m}dz_m\Big)^{\frac{r_{m-1}'}{r_{m}'}}\\
&&\qquad
\dots|z_2|^{r_{2}'s_{2}}dz_2\Big)^{\frac{r_1'}{r_2'}}|z_1|^{r_1's_1}dz_1\Big)^{\frac{1}{r_1'}}
2^{\kappa}R(2^{j_0}R)^{-s_1} \prod_{i=2}^m(2^{j_i}R)^{-s_i}2^{\kappa mn}\\
&\lesssim
&\sum_{|\alpha|=1}\Big(\int_{C_{j_0}}\Big(\int_{S_{j_2}(B_R)}\dots\Big(\int_{S_{j_m}(B_R)}
\Big|{\cal F}^{-1}(\xi_1^\alpha{\sigma}_{\kappa}) (
2^{\kappa}z_1,\dots,2^{\kappa}z_m)\Big|^{r_m'}| z_m|^{r_m's_m}dz_m\Big)^{\frac{r_{m-1}'}{r_{m}'}}\\
&&\qquad \dots| z_{2}|^{r_2's_2}dz_2\Big)^{\frac{r_1'}{r_2'}}|
z_1|^{r_1's_1}dz_1\Big)^{\frac{1}{r_1'}}
2^{\kappa}R(2^{j_0}R)^{-s_1}
\prod_{i=2}^m(2^{j_i}R)^{-s_i}2^{\kappa mn}\\
&\lesssim
&\sum_{|\alpha|=1}\|\xi_1^{\alpha}\sigma_{\kappa}\|_{W^{s_1,\dots,s_m}(\mathbb{R}^{mn})}2^{\kappa}R(2^{j_0}R)^{-s_1}
\prod_{i=2}^m(2^{j_i}R)^{-s_i}2^{-\kappa(s_1+\dots+s_m-n/r_1-\dots-n/r_m)}\\
&\lesssim &2^{\kappa}R(2^{j_0}R)^{-s_1}
\prod_{i=2}^m(2^{j_i}R)^{-s_i}2^{-\kappa(s_1+\dots+s_m-n/r_1-\dots-n/r_m)},
\end{eqnarray*}
where $C_{j_0}=\{z:\, 2^{j_0-2}R\leq |z|\leq 2^{j_0+1}R\}.$

Step 2: if $\min\{j_0,\,j_2,\dots,j_m\}=0$, for example,
$j_{l+1}=\dots=j_m=0$ and $j_{k}\geq 1$ for $2\leq k\leq l$, then, for index $\alpha\in \mathbb{Z}_+^n$ and $2^{\kappa} R<1$,
\begin{eqnarray*}
&&\Big(\int_{C_{j_0}}\Big
(\int_{S_{j_2}(B_R)}\dots\Big(\int_{S_{j_m}(B_R)}
\Big|\partial^{\alpha,\,0,\dots,0}{\cal F}^{-1}{\sigma}_{\kappa} (
2^{\kappa}z_1,\dots,2^{\kappa}z_m)\Big|^{r_m'}dz_m\Big)^{\frac{r_{m-1}'}{r_{m}'}}dz_{m-1}\\
&&\qquad \dots\Big)^{\frac{r_1'}{r_2'}}dz_1\Big)^{\frac{1}{r_1'}}2^{\kappa nm}2^{\kappa}R\\
&&\qquad\lesssim
\Big(\int_{C_{j_0}}\Big(\int_{S_{j_2}(B_R)}\dots\Big(\int_{S_{j_l}(B_R)}\Big
(\int_{B_R}\dots\Big(\int_{B_R} \Big|{\cal
F}^{-1}(\xi_1^{\alpha}\sigma_{\kappa}) (
2^{\kappa}z_1,\dots,2^{\kappa}z_m)\Big|^{r_m'}dz_m        \Big)^{\frac{r_{m-1}'}{r_{m'}}}  \\
&&\qquad\qquad \dots
\Big)^{\frac{r_{l}'}{r_{l+1}'}}|z_l|^{r_l's_l}dz_l\Big)\dots\Big)^{\frac{r_1'}{r_2'}}|
z_1|^{r_1's_1}dz_1\Big)^{\frac{1}{r_1'}}
2^{\kappa}R(2^{j_0}R)^{-s_1} \prod_{i=2}^l(2^{j_i}R)^{-s_i}2^{\kappa mn}\\
&&\qquad\lesssim 2^{\kappa}R(2^{j_0}R)^{-s_1}
\dprod_{i=2}^l(2^{j_i}R)^{-s_i}2^{-\kappa(s_1+\dots+s_l-n/r_1-\dots-n/r_m)}\\
&&\qquad\lesssim 2^{\kappa}R (2^{j_0}R)^{-s_1}
\dprod_{i=2}^l(2^{j_i}R)^{-s_i}R^{-s_{l+1}-\dots-s_m}2^{-\kappa(s_1+\dots+s_m-n/r_1-\dots-n/r_m)}.
\end{eqnarray*}
This together with the argument in step 1 shows that, for $j_2,\dots,j_{l}$ positive and
$j_{l+1}=\dots=j_m=0$, inequality (\ref{eq3.1}) holds $k=1$.
\qed

\begin{lem}\label{lem3.4}
Let $\sigma$ be a multiplier satisfying (\ref{eq1.7}) for $s_1,\dots,s_m\in (n/2,\,n]$, and $r_1,\dots,r_m\in (1,\,2]$ with $r_l>n/s_l$ for $1\leq l\leq m$. Then, for each $\kappa\in \mathbb{Z}$ and integer $k$ with $1\leq k\leq m$, there exists a function ${\rm H}_{k}^{\kappa}$ such that, for any ball $B$ with radial $R$, any function $f_k$ with ${\rm supp}\,f_k\subset B$, and any $x\in \rn\backslash 4B$, $y_k,y_k'\in B$,
\begin{eqnarray}\label{eq3.2}
\begin{array}[b]{cl}&\dint_{(\rn)^m}|W_{k,\,\kappa}(x,y_1,\dots,y_m;\,y_k)|\prod_{l=1}^m|f_l(y_l)|d\vec{y}\\
&\qquad\lesssim \dint_{\rn}|f_k(y_k)|{\rm H}_{k}^{\kappa}(x,\,y_k,\,y_k')dy_k \prod_{1\leq l\leq m,\,l\not =k}M_{r_l}f_l(x)
\end{array}
\end{eqnarray}
and
\begin{eqnarray}\label{eq3.3}
\Big(\int_{S_{j_0}(B)}|{\rm
H}^{\kappa}_{k}(x,\,y_k,\,y_k')|^{r'_k}dx\Big)^{1/r'_k}\lesssim
\frac{2^{-\kappa(s_k-n/r_k)}}{(2^{j_0}R)^{s_k}}\qquad\text{for
integer}\ j_0\geq 3.
\end{eqnarray}
\end{lem}

\noindent
{\it Proof}. We show the case $k=1$ only. Let
\begin{eqnarray*}
{\rm H}_1^{\kappa,\,1}(x,\,y_1) &=&2^{\kappa
n}\Big(\int_{\rn}\Big(\int_{\rn}\dots\Big(\int_{\rn}|{\cal
F}^{-1}\sigma_{\kappa}(2^{\kappa}(x-y_1),2^{\kappa}x-y_2,\dots,2^{\kappa}x-y_m)|^{r_m'}\\
&&\qquad\times\langle
2^{\kappa}x-y_m\rangle^{r_m's_m}dy_m\Big)^{\frac{r_{m-1}'}{r_m'}}\dots\Big)^{\frac{r_2'}{r_3'}}\langle
2^{\kappa}x-y_2\rangle^{r_2's_2}dy_2\Big)^{\frac{1}{r_2'}}.
\end{eqnarray*}
For any integer $j_0\geq 3$ and $y_1\in B$, we deduce from Lemma \ref{lem3.1} that
\begin{eqnarray*}
&&\Big(\int_{S_{j_0}(B)}|{\rm
H}_1^{\kappa,\,1}(x,\,y_1)|^{r_1'}dx\Big)^{\frac{1}{r_1'}}\\
&&\qquad=2^{\kappa
n}\Big(\int_{S_{j_0}(B)}\Big(\int_{\rn}\Big(\dots\Big(\int_{\rn}|{\cal
F}^{-1}\sigma_{\kappa}(2^{\kappa}(x-y_1),2^{\kappa}x-y_2,\dots,2^{\kappa}x-y_m)|^{r_m'}\\
&&\qquad\quad\times\langle
2^{\kappa}x-y_m\rangle^{r_m's_m}dy_m\Big)^{\frac{r_{m-1}'}{r_m'}}\dots\Big)^{\frac{r_2'}{r_3'}}\langle
2^{\kappa}x-y_2\rangle^{r_2's_2}dy_2\Big)^{\frac{r_1'}{r_2'}}\\
&&\qquad\quad\times|2^{\kappa
}(x-y_1)|^{r_1's_1}dx\Big)^{\frac{1}{r_1'}}(2^{\kappa}2^{j_0}R)^{-s_1}\\
&&\qquad\lesssim 2^{\kappa
n}\Big(\int_{\rn}\Big(\int_{\rn}\Big(\dots\Big(\int_{\rn}|{\cal
F}^{-1}\sigma_{\kappa}(z_1,\dots,z_m)|^{r_m'}\langle z_m\rangle^{r_m's_m}dz_m\Big)^{\frac{r_{m-1}'}{r_m'}}\dots\Big)^{\frac{r_2'}{r_3'}}\\
&&\qquad\quad \times\langle
z_2\rangle^{r_2's_2}dz_2\Big)^{\frac{r_1'}{r_2'}}\langle
z_1\rangle^{r_1's_1}dz_1\Big)^{\frac{1}{r_1'}}2^{-\kappa
n/r_1'}(2^{\kappa}2^{j_0}R)^{-s_1}\\
&&\qquad\lesssim2^{-\kappa (s_1-n/r_1)}(2^{j_0}R)^{-s_1}.
\end{eqnarray*}
Since $r_ls_l>n$, H\"older's inequality leads to
\begin{eqnarray*}
&&\int_{(\rn)^{m}}|{\cal
F}^{-1}\widetilde{\sigma}_{\kappa}(x-y_1,\dots,x-y_m)||f_1(y_1)\dots f_m(y_m)|d\vec{y}\\
&&\qquad=2^{\kappa mn}\int_{(\rn)^{m}}|{\cal
F}^{-1}\sigma_{\kappa}(2^{\kappa}(x-y_1),\dots,2^{\kappa}(x-y_m))||f_1(y_1)\dots f_m(y_m)|d\vec{y}\\
&&\qquad\lesssim 2^{\kappa
mn}\int_{\rn}\Big(\int_{\rn}\Big(\int_{\rn}\dots\Big(\int_{\rn}|{\cal
F}^{-1}\sigma_{\kappa}(2^{\kappa}(x-y_1),\dots,2^{\kappa}(x-y_m))|^{r_m'}\\
&&\qquad\quad\times\langle 2^{\kappa}(x-y_m)\rangle^{r_m's_m}d
y_m\Big)^{\frac{r_{m-1}'}{r_m'}}\dots\Big)^{\frac{r_2'}{r_3'}}\langle
2^{\kappa}(x-y_2)\rangle^{r_2's_2}d
y_2\Big)^{\frac{1}{r_2'}}|f_1(y_1)|dy_1\\
&&\qquad\quad\times\prod_{l=2}^m\Big(\int_{\rn}\frac{|f_l(y_l)|^{r_l}}{\langle2^{\kappa}(x-y_l)\rangle^{r_ls_l}}dy_l\Big)^{1/r_l}\\
&&\qquad\lesssim\int_{\rn}{\rm
H}^{\kappa,\,1}_1(x,\,y_1)|f_1(y_1)|\,dy_1\prod_{l=2}^mM_{r_l}f_l(x).
\end{eqnarray*}
We can also verify that for $y_1'\in B$,\begin{eqnarray*}
&&\int_{(\rn)^{m}}|{\cal
F}^{-1}\widetilde{\sigma}_{\kappa}(x-y_1',\dots,x-y_m)||f_1(y_1)\dots f_m(y_m)|d\vec{y}\\
&&\qquad\lesssim{\rm H}_1^{\kappa,
1}(x,\,y_1')\int_{\rn}|f_1(y_1)|\,dy_1\prod_{l=2}^mM_{r_l}f_l(x)
\end{eqnarray*} and
\begin{eqnarray*}
\Big(\int_{S_{j_0}(B)}|{\rm
H}^{\kappa,\,1}_{1}(x,\,y_1')|^{r_1}dx\Big)^{\frac{1}{r_1}}\lesssim
\frac{2^{-\kappa(s_k-n/r_k)}}{(2^{j_0}R)^{s_k}}\qquad\text{for
integer}\ j_0\geq 3.
\end{eqnarray*}
Taking ${\rm H}_{1}^{\kappa}(x,\,y_1,\,y_1')={\rm
H}_1^{\kappa,\,1}(x,\,y_1)+{\rm H}_1^{\kappa,\,1}(x,\,y'_1)$, we
complete the proof. \qed

\begin{lem}\label{lem3.5}
Let $m$ and $k$ be positive integers with $k\leq m$, $\sigma$ be a multiplier satisfying (\ref{eq1.7}) for some $s_1,\dots,s_m\in (n/2,\,n]$, and $r_1,\dots,r_m\in (1,\,2]$.
Then, for any ball $B$ with radial $R$, $x,x'\in \frac{1}{4}B$, and nonnegative integers $j_1,\dots,j_m$,
$$\begin{array}[b]{cl}
&\Big(\dint_{S_{j_1}(B)}\Big(\int_{S_{j_2}(B)}\dots\Big(\dint_{S_{j_m}(B)}|W_{0,\,\kappa}
(x,\,y_1,\dots,y_m;\,x')|^{r_m'}dy_m\Big)^{\frac{r_{m-1}'}{r_m'}}\dots\Big)^{\frac{r_1'}{r_2'}}dy_1\Big)^{\frac{1}{r_1'}}\\
&\qquad\lesssim R
\dfrac{2^{-\kappa(s_1+\dots+s_m-n/r_1-\dots-n/r_m-1)}}{\prod_{1\leq
i\leq m}(2^{j_i}R)^{s_i}}
\end{array}$$
provided $2^{\kappa}R<1$.
\end{lem}

This lemma can be obtained by the argument used in the proof of Lemma \ref{lem3.3}.

\begin{lem}\label{lem3.6}
Let $m$ and $k$ be positive integers with $k\leq m$, and $\sigma$ be a multiplier satisfying (\ref{eq1.7}) for some $s_1,\dots,s_m\in (n/2,\,n]$. Suppose that $r_1,\dots,r_m\in (1,\,2]$ such that $r_k>n/s_k$ for
$k=1,\dots,m$. Then, for any $B$ with radial $R$, $x,x'\in B$, integer $j_k\geq 2$, and functions $f_1,\dots,f_m$ satisfying ${\rm supp}\,f_k\subset \rn\backslash 4B$ for some $k\in \{1,\dots, m\}$,
\begin{eqnarray*}
&&\int_{S_{j_k}(B)}\Big(\int_{(\rn)^{m-1}}|W_{0,\,\kappa}(x,\,y_1,\dots,y_m;
x')|\prod_{l=1}^m|f_l(y_l)|dy_1\dots dy_{k-1}dy_{k+1}dy_m\Big)dy_k\\
&&\qquad\lesssim\frac{2^{-\kappa(s_k-n/r_k)}}{(2^{j_k}R)^{s_k-n/r_k}}\prod_{l=1}^{m}\Big(M_{r_l}f_l(x)+M_{r_l}f_l(x')\Big).
\end{eqnarray*}
\end{lem}

\noindent
{\it Proof}. We consider the case $k=1$ only. As in the proof of Lemma \ref{lem3.4}, we have
\begin{eqnarray*}
&&\int_{S_{j_1}(B)}\Big(\int_{(\rn)^{m-1}}|{\cal
F}^{-1}\widetilde{\sigma}_{\kappa}(x-y_1,\dots,x-y_m)|\prod_{l=2}^m|f_l(y_l)|dy_2\dots dy_m\Big)|f_1(y_1)|dy_1\\
&&\qquad\lesssim\frac{2^{-\kappa(s_1-n/r_1)}}{(2^{j_1}R)^{s_1-n/r_1}}\prod_{l=1}^{m}M_{r_l}f_l(x)
\end{eqnarray*}
and
\begin{eqnarray*}
&&\int_{S_{j_1}(B)}\Big(\int_{(\rn)^{m-1}}|{\cal
F}^{-1}\widetilde{\sigma}_{\kappa}(x'-y_1,\dots,x'-y_m)|\prod_{l=2}^m|f_l(y_l)|dy_2\dots dy_m\Big)|f_1(y_1)|dy_1\\
&&\qquad\lesssim\frac{2^{-\kappa(s_1-n/r_1)}}{(2^{j_1}R)^{s_1-n/r_1}}\prod_{l=1}^{m}M_{r_l}f_l(x').
\end{eqnarray*}
So we get the desired conclusion directly.
\qed
\medskip

We return to show Theorem \ref{thm1.2}.

\medskip
\noindent
{\it Proof of Theorem \ref{thm1.2}}.
We will employ the argument given in \cite[p.\,350]{kw}. For $N\in\mathbb N$, let
$$\sigma^N(\xi_1,\dots,\xi_m)=\sum_{|\kappa|\leq
N}\widetilde{\sigma}_{\kappa}(\xi_1,\dots,\xi_m)$$ and denote by
$T_{\sigma,\,N}$ the multiplier operator associated with $\sigma^N$.
It is obvious that $T_{\sigma,\,N}$ is an $m$-linear singular
integral operator with kernel
$$K^N(x;\,y_1,\dots,y_m)={\cal F}^{-1}\sigma^N(x-y_1,\dots,x-y_m)$$ in the sense of (\ref{eq1.1}). Note that for
$f_1,\dots,f_m\in \mathscr{S}(\rn)$,
$$\|T_{\sigma}(f_1,\dots,f_m)-T_{\sigma,\,N}(f_1,\dots,f_m)\|_{L^{\infty}(\rn)}\lesssim
\|(\sigma-\sigma^{N})\widehat{f_1}\dots\widehat{f_m}\|_{L^1(\rn)}\to 0.$$
as $N\to\infty$. By a density argument, it
suffices to prove that the conclusions of Theorem \ref{thm1.2} are true for
$T_{\sigma,\, N}$ with bound independent of $N$.

Let $t_k=n/s_k$. We only need to show that, when $\sigma$ satisfies
(\ref{eq1.7}) for $s_1,\dots,s_m\in (n/2,\,n]$, all of the
assumptions (i)$-$(iv) in Theorem \ref{thm1.1} hold for the operator
$T_{\sigma,\,N}$ provided, for $k=1,\dots, m$, each $r_k\in (t_k,\,2)$ closes enough to $t_k$ (satisfy $n/r_k>s_k-1/m$).
By Lemma \ref{lem3.1}, for $x\in \rn$ and integers $j_1,\,\dots \,j_m\in\{1,\,2\}$,
\begin{eqnarray*}
&&\Big(\dint_{S_{j_1}(B(x,\,R))}\Big(\dots\Big(\dint_{S_{j_m}(B(x,\,R))}|{\cal F}^{-1}\widetilde{\sigma}_{\kappa}(x-y_1,\dots,x-y_m)|^{r_m'}dy_m\Big)^{\frac{r_{m-1}'}{r_m'}}\dots\Big)^{\frac{r_1'}{r_2'}}dy_1\Big)^{\frac{1}{r_1'}}\\
&&\qquad\lesssim
\Big(\int_{S_{j_1}(B_R)}\Big(\dots\Big(\int_{S_{j_m}(B_R)}|{\cal
F}^{-1}\widetilde{\sigma}_{\kappa}(z_1,\dots,z_m)|^{r_m'}\langle z_m\rangle^{r_m'\frac{n}{r_m}}dz_m\Big)^{\frac{r_{m-1}'}{r_m'}}\dots\Big)^{\frac{r_1'}{r_2'}}\\
&&\qquad\qquad\times
\langle z_1\rangle^{r_1'\frac{n}{r_1}}dz_1\Big)^{\frac{1}{r_1'}}R^{-n/r_1-\dots-n/r_m}\\
&&\qquad\lesssim 2^{\kappa n(1/r_1+\dots+1/r_m)}\|\sigma_{\kappa}\|_{W^{n/r_1,\dots
n/r_m}(\mathbb{R}^{mn})}R^{-n/r_1-\dots-n/r_m},
\end{eqnarray*}
which implies
\begin{eqnarray*}&&\Big(\dint_{A_R^x}\Big(\int_{A_R^x}\dots\Big(\dint_{A_R^x}|K^N(x;y_1,\dots,y_m)|^{r_m'}dy_m\Big)^{\frac{r_{m-1}'}{r_m'}}\dots\Big)^{\frac{r_1'}{r_2'}}dy_1\Big)^{\frac{1}{r_1'}}\\
&&\qquad\lesssim 2^{n N(1/r_1+\dots+1/r_m)}R^{-n/r_1-\dots-n/r_m},
\end{eqnarray*}
and hence $T_{\sigma,\,N}$ satisfies assumption (i) of Theorem \ref{thm1.1}.
Denote
$$W_{k}^N(x,\,y_1,\dots,y_m;\,y_k')=K^N(x;\,y_1,\dots,y_m)-K^N(x;\,y_1,\dots,y_k',\dots,\,y_m)$$
for $k=1,\dots,m$, and
$$W_{0}^N(x,\,y_1,\dots,y_m;\,x')=K^N(x;\,y_1,\dots,y_m)-K^N(x';\,y_1,\dots,y_m).$$
Let $B$ be a ball with radial $R$. For $x,\,x'\in B$ and for functions $f_1,\dots,f_m$
with ${\rm supp}\,f_k\subset \rn\backslash 4B$, it follows from Lemmas \ref{lem3.5} and \ref{lem3.6} that
\begin{eqnarray*}
&&\int_{(\rn)^m}|W_0^N(x,y_1,\dots,\,y_m,\,x')||f_1(y_1)\dots f_m(y_m)|d\vec{y}\\
&&\qquad\lesssim\sum_{\kappa:\,2^{\kappa}R>1}\sum_{j_k=2}^{\infty}\int_{S_{j_k}(B)}
\int_{(\rn)^{m-1}}|W_{0,\,\kappa}(x,y_1,\dots,\,y_m,\,x')||f_1(y_1)\dots f_m(y_m)|d\vec{y}\\
&&\qquad\quad+\sum_{\kappa:\,2^{\kappa}R\leq
1}\sum_{j_1,\dots,j_m=0}^{\infty}\Big(\dint_{S_{j_1}(B)}\Big(\dots\Big(\dint_{S_{j_m}(B)}|W_{0,\,\kappa}
(x,\,y_1,\dots,y_m;\,x')|^{r_m'}dy_m\Big)^{\frac{r_{m-1}'}{r_m'}}\\
&&\qquad\qquad \dots\Big)^{\frac{r_1'}{r_2'}}dy_1\Big)^{\frac{1}{r_1'}}\prod_{k=1}^mM_{r_l}f_l(x)\prod_{l=1}^m|S_{j_l}(B)|^{\frac{1}{r_l}}\\
&&\qquad\lesssim\sum_{\kappa:\,2^{\kappa}R>1}(2^{\kappa}R)^{-s_k+n/r_k}\prod_{l=1}^m\Big(M_{r_l}f_l(x)+M_{r_l}f_l(x')\Big)\\
&&\qquad\quad +\sum_{\kappa:\,2^{\kappa}R\leq
1}(2^{\kappa}R)^{-s_1-\dots-s_m+n/r_1+\dots+n/r_m+1}\prod_{l=1}^mM_{r_l}f_l(x)\\
&&\qquad\lesssim\prod_{l=1}^m\Big(M_{r_l}f_l(x)+M_{r_l}f_l(x')\Big),
\end{eqnarray*}
since $s_1+\dots+s_m<n/r_1+\dots+n/r_m+1$. Thus, $T_{\sigma,\,N}$ satisfies assumption (ii) of Theorem \ref{thm1.1}. For a fixed integer $k\in \{1,\dots, m\}$, a ball $B$ with radial $R$, and $x,\,y_k,\,y_k'\in \rn$, set
\begin{eqnarray*}
{\rm H}_{k,\,B}(x,\,y_k,\,y_k')&=&\sum_{\kappa:\,2^\kappa
R<1}\Big(\int_{\rn}\dots\Big(\int_{\rn}|W_{k,\,\kappa}
(x,y_1,\dots,y_m;y_k')|^{r_m'}dy_m\Big)^{\frac{r_{m-1}'}{r_m'}}\dots dy_1\Big)^{\frac{1}{r_1'}}\\
&&\quad+\sum_{\kappa:\, 2^\kappa R>1}{\rm H}_{k}^{\kappa}(x,\,y_k,\,y_k'),
\end{eqnarray*}
where ${\rm H}_{k}^{\kappa}$ is the function satisfying $(\ref{eq3.2})-(\ref{eq3.3})$ and the integral on the right hand side of equality is taken with respect to variables $y_1,\dots, y_{k-1}, y_{k+1},\dots, y_1$.
By Lemmas \ref{lem3.3} and \ref{lem3.4}, we see that $T_{\sigma,\,N}$ satisfies assumption (iii) of Theorem \ref{thm1.1}. Finally, it follows from
\cite[Theorem 1.1]{gmt} that $T_{\sigma,\,N}$ is bounded from $L^2(\rn)\times
L^{\infty}(\rn)\times\dots\times L^{\infty}(\rn)$ to $L^2(\rn)$ provided $s_1,\dots,s_m>n/2$, and hence assumption (iv) holds. This finishes the proof of Theorem \ref{thm1.2}.
\qed
\medskip

\noindent
{\bf Acknowledgement.} The authors would like to thank Professor N.
Tomita for helpful discussion and giving preprints \cite{ft,gmt,mt} to the authors.

\end{document}